\documentclass[final,1p,times,authoryear]{elsarticle}
\usepackage{amsmath}
\usepackage{amsthm}
\usepackage{amssymb}
\usepackage{amsfonts}
\usepackage{natbib}
\setcitestyle{numbers,square}

\begin{document}

\begin{frontmatter}

\title{Automated Proof of Mixed Trigonometric-polynomial Inequalities
in the Unbounded Case}

\author{Chen shiping}
\address{Sichuan Trade School, Yaan, 625107, Sichuan, China}
\ead{chinshiping@sina.com}
\ead[url]{URL 1}

\author{Ge Xinyu}
\address{Chengdu Institute of Computer Application, Chinese Academy of Sciences, Chengdu, 610000, Sichuan, China}
\address{University of Chinese Academy of Sciences, Beijing, 100049, Beijing, China}
\ead{geeexy@163.com}
\ead[url]{URL 2}

\begin{abstract}
Mixed trigonometric-polynomials frequently occur in applications in physics, numerical analysis and engineering, the algorithm has been already proposed to determine its sign on \ $(0,\frac{\pi}{2}]$. This paper proposes a procedure to extend the interval to \ $(0,+\infty)$. Such generalization is by no means trivial, for the process depends on boundedness of roots and square-free factorization of mixed trigonometric-polynomials.
\end{abstract}
\begin{keyword}
Mixed trigonometric-polynomial; Boundedness of roots; Square-free factorization; Automated proof of inequalities.
\end{keyword}
\end{frontmatter}

\section{Introduction}
Automatic proof of the transcendental inequalities and sign-deciding of transcendental functions have always been hot and difficult problems in the field of mathematical mechanization and automatic reasoning in recent years. Among the transcendental functions, a class of functions named mixed trigonometric-polynomial (denoted as Trigonometric-polynomial or MTP simply) and defined by the formula \ $F(x)=\sum a_i x^{p_i} \sin^{q_i}(\omega_i x)\cos^{r_i}(\omega_i x)$ attract more and more scholars' attention$^{[1-8]}$, which frequently occur in applications in physics, numerical analysis and engineering, where \ $a_i\in \mathbb{R}$,\ $p_i,q_i,r_i\in \mathbb{N}\cup \{0\}$,\ $\omega_i\in \mathbb{Q}$. Moreover, an inequality of the form \ $F(x)>0$ is called mixed trigonometric polynomial inequality (MTP inequality), where \ $F(x)$ is an MTP. Let \ $p$ be the least common divisor of numerators of \ $\omega_i$s, and \ $q$ the greatest common multiple of denominators, \ $y=\frac{p}{q} x$, then an MTP can be reduced to another equivalent form \ $G(y)=g(y,\sin (y),\cos(y))= \sum a_i y^{p_i} \sin^{q_i}(y)\cos^{r_i}(y)$, where \ $g\in \mathbb{R}[x_1,x_2,x_3]$.

Ref[8] presents an algorithm for automated proof of MTP inequality and deciding the sign of MTP. By Taylor expansion, the proving of the target inequality is reduced to the verification of a series of polynomial inequalities with only one variable, and then completed by algebraic inequality-proving package such as BOTTEMA. The algorithm is complete and experiments show that it is very effective for mixed trigonometric-polynomial inequalities, furthermore, the procedure is 'readable'. However, the algorithm limits the domain of inequality on \ $(0,\frac{\pi}{2}]$ because the Taylor expansion of \ $\arctan(x)$ is convergent only on \ $[0,1]$, where it is assumed that \ $x\geq0$. Of course, the domain of the trigonometric-polynomial inequalities can be extended to\ $(0,T)$ for \ $\forall T>0$ by the double angle formulas of trigonometric function.

In this paper, we are to discuss how to decide the sign of mixed trigonometric-polynomial \ $f(x) =f(x,\sin(x),\cos(x))$ on the unbounded interval \ $(0,+\infty)$. The generalization is non-trivial, for it is necessary to discuss the boundedness of MTP's roots and the square-free decomposition of  MTPs to deal with the unbounded case.

Ref[9] discusses the boundedness of roots of the so-called trigonometric-exponential polynomial with the form of  \ $f(t)=\sum e^{r_j t}(P_{1,j}(t)\cos (\omega_j t) + P_{2,j} (t) \sin (\omega_j t))$, and presents a deciding procedure in two cases using the theory of semi-algebraic sets , the first case is that the set of frequencies ${\omega_j}$ spans a one-dimensional vector space over \ $\mathbb{Q}$ , the second is that  the set of frequencies ${\omega_j}$ spans a two-dimensional vector space over \ $\mathbb{Q}$ and the polynomials \ $P_{1,j}(t)$ and $P_{2,j}(t)$ are all constants. Ref[9]'s methods can be used to decide the boundedness of MTP's roots, while we will design a simpler and more convenient scheme for the specific MTPs with the help of the Sylvester resultant.

Factorization of trigonometric functions is a classic field and the following methods are frequently used. The first is the quotient ring \ $Q[s,c]/<s^2+c^2-1>$, but it is not a unique factorization domain and so, the factorization is not unique in general,furthermore,it is still needed to decide whether each factor has multiple roots$^{[10]}$. The second is using Tan-half angle substitution, \ $\sin (x)=\frac{2t}{1+t^2}$, $\cos (x)=\frac{1-t^2}{1+t^2}$, to transform the trigonometric polynomial to a rational expression, where $t=\tan(\frac{x}{2})$ and \ $x\neq (2k+1)\pi$ for \ $k\in \mathbb{N}$ $^{[10-11]}$. The third is reducing the trigonometric functions to polynomials in complex field \ $Q(I)[e,e^{-1}]$ by Euler Theorem$^{[11]}$, where \ $e$ denotes \ $e^{Ix}$. This scheme requires \ $g(e)\in Q(I)[e,e^{-1}]$ be even or odd function of \ $e$. Besides, the above methods can deal with the trigonometric functions with form \ $f(\sin (x),\cos (x))$ only, where \ $f(x,y)\in \mathbb{R}[x,y]$ or \ $\mathbb{Q}[x,y]$. In this paper, we will present a procedure to decompose the mixed trigonometric-polynomials with the form \ $f(x,\sin (x),\cos (x))$, i.e. the monomials contain one variable and trigonometric functions applied to the same variable, and the domain is \ $(0,+\infty)$ or \ $(-\infty,+\infty)$ without excluding any special point.

The rest of the paper is organized as follows. Section 2 studies the positive root boundedness of MTP. Section 3 proposes the scheme of square-free factorization of MTP. Section 4 presents deciding procedure of the sign of MTP on unbounded interval $(0,+\infty)$. We will conclude the paper in Section 5.

\section{Decision of the Positive Root Boundedness of Mixed Trigonometric-Polynomial}

In this section, we aim to decide the positive root boundedness of mixed trigonometric-polynomial \ $F(x)=f(x,\sin (x),\cos (x))$.

\noindent\textbf{Algorithm 2.1} Decision of the positive root boundedness of MTP

Input: \ $F(x)=f(x,\sin (x),\cos (x))$, where \ $f\in \mathbb{R}_{alg}[x,y,z]$; Where \ $\mathbb{R}_{alg}$ is the set of real algebraic numbers.

Output: \ whether the positive roots of \ $F(x)$ are bounded or not, if bounded, output one upper bound.

1)\ $f\_{lst}\gets factor(square-free(f(x,s,c)))$;

2)\ To decide if positive roots of each function in \ $f\_{lst}$ are bounded, \ $f_i(x,s,c)\gets f\_{lst}[i]$; Write \ $f_i(x,s,c)$ as $f(x,s,c)$ for convenience below.

\ \ \ 2.1)\ if \ $degree(f,c)$ is odd, then

\ \ \ \ \ 2.1.1)$f_1(x,s)\gets square-free(res(f,s^2+c^2-1,c))$,

\ \ \ \ \ 2.1.2)\ $f_2(x)\gets square-free(res(f_1,diff(f_1,s),s))$;

\ \ \ \ \ 2.1.3)\ $r_0\gets$ the maximum positive root of \ $f_2(x)$. (\ $r_0\gets 0$ if \ $f_2(x)$ has no positive root );

\ \ \ \ \ 2.1.4)\ Select \ $\forall r_1>r_0$, if \ $g(s)=f_1(r_1,s)=0$ has no real root, then the positive root of \ $f(x,\sin(x),\cos(x))$ is bounded and \ $r_0$ is an upper bound.

\ \ \ \ \ \ \ If \ $g(s)=f_1(r_1,s)=0$ has at least one real root, then the positive root of \ $f(x,\sin(x),\cos(x))$ is unbounded, which implies that the positive root of \ $F(x)$ is unbounded. The algorithm ends.

\ \ \ 2.2)Else  (the case that $degree(f,c)$ is even)

\ \ \ \ \ 2.2.1)\ Reduce \ $f(x,s,c)$ to \ $g(x,s)$ by the substitution \ $c^2=1-s^2$; where we suppose that \ $degree(g,s)>0$, otherwise \ $g$ is a unary polynomial of \ $x$ and its positive roots are bounded.

\ \ \ \ \ 2.2.2)\ $g_1(x)\gets res(s+1,g,s)\times res(s-1,g,s)\times res(g,diff(g,s),s)$;

\ \ \ \ \ 2.2.3)\ $r_0\gets$ the maximum positive root of \ $g_1(x)$. (\ $r_0\gets 0$ if \ $g_1(x)$ has no positive root );

\ \ \ \ \ 2.2.4)\ Select \ $\forall r_1>r_0$, if \ $g_2(s)=g(r_1,s)$ has no real root on \ $(-1,1)$, then the positive root of \ $g(x,\sin(x))$ is bounded,  and \ $r_0$ is an upper bound.

\ \ \ \ \ \ \ If  \ $g_2(s)=g(r_1,s)$ has a least one real root on \ $(-1,1)$, then the positive root of \ $g(x,\sin(x))$ is unbounded, i.e. the positive root of \ $F(x)$ is unbounded. The algorithm ends.

3)\ The maximum value of the positive roots bounds of all functions in \ $f\_lst$ is one bound of \ $F(x)$. The algorithm ends.

Here \ $res(f,g,x)$ denotes the Sylvester resultant of \ $f(x)$ and \ $g(x)$ about \ $x$, similarly hereinafter.

\noindent\textbf{Lemma 2.1}$^{[12]}$
Suppose that \ $f(a,x)=a_mx^m+a_{m-1}x^{m-1}+\cdots+a_0$ is a unary polynomial of real parameter  coefficients, where \ $a$ represents \ $a_m, a_{m-1},\dots, a_0$. Denote \ $R(a)=res(f,f',x)$, if \ $c_1$ and \ $c_2$ are two points in the same connected branch of \ $R(a)\neq 0$ in parameter space \ $R^{m+1}$, then \ $f(c_1,x)$ and \ $f(c_2,x)$ have the same number of real roots.

\noindent\textbf{Theorem 2.1} Algorithm 2.1 is correct.
\begin{proof}
First, suppose that \ $degree(f,c)$ is odd.

If \ $\exists r_1>r_0$,\ $g(s)=f_1(r_1,s)=0$ has no real root, then \ $g_r(s)=f_1(r,s)=0$ has no real root for any \ $r>r_0$ by Lemma 2.1. Specially, \ $\sin(r)$ is not a root of \ $g_r(s)$, so \ $f_1(r,\sin(r))\neq 0$, which is to say that \ $res(f(r,\sin (r),c),\sin^2(r)+c^2-1,c)\neq 0$, so \ $f(r,\sin(r),c)$ has no common root with \ $\sin^2(r)+c^2-1$, as \ $\sin^2(r)+\cos^2(r)-1=0$, so \ $f(r,\sin (r),\cos (r))\neq 0$, hence the root of \ $f(x,\sin (x),\cos (x))$ is not bigger than \ $r_0$, so its positive roots is bounded.

If $\exists r_1>r_0$,\ $g(s)=f_1(r_1,s)=0$ has at least one real root. Suppose \ $T$ be an arbitrary given positive real number.

Let \ $c=0$ under the assumption that \ $f(x,s,c)=0$ and \ $s^2+c^2=1$, \ $s=1$ or \ $-1$ obviously. If \ $f(x,1,0)$ or \ $f(x,-1,0)$ is identically equal to \ $0$, the positive roots of \ $f(x,\sin (x),\cos (x))$ is unbounded clearly. Without loss of generality, we can suppose that \ $f(x,1,0)$ and \ $f(x,-1,0)$ are non-constant polynomials of \ $x$, their biggest positive roots are \ $r_1$ and \ $r_2$ respectively, then we have that \ $c\neq 0$ for \ $\forall x>r_3=max\{r_1,r_2\}$, i.e.\ $c$ has a constant sign for \ $x>r_3$.

By Lemma 2.1, for any \ $r>T_0=max\{r_3,T,r_0\}$, where \ $r_0$ is the maximum positive root of \ $f_2(x)$ in Algorithm 2.1, \ $g_r(s)=f_1(r,s)=0$ has at least one real root, denotes \ $s(r)$ as its smallest one, then by the continuous dependence of the roots of a polynomial on its coefficients, we know that \ $s(r)$ is continuous on \ $(T_0,+\infty)$.

Recall $f_1(r,s(r)) = 0$ implies that \ $f(r,s(r),c)$ and \ $s(r)^2+c^2-1$ have common roots, denoted by \ $c_0$, by the assumption that\ $degree(f,c)$ is odd we have that \ $c_0$ must be real-valued, so we can declare that \ $-1\leq s(r)\leq 1$ due to \ $s(r)^2+c_0^2-1=0$.

Assume that \ $c>0$ for \ $x\in(r_3,+\infty)$. For \ $\forall k\in \mathbb{N}$, the value of \ $\sin (x)$ changes from \ $-1$ to \ $1$ continuously on \ $[2k\pi-\frac{\pi}{2},2k\pi+\frac{\pi}{2}]$. Suppose \ $k$ is big enough such that \ $2k\pi-\frac{\pi}{2}>T_0$, then \ $\sin (x)$ and \ $s(r)$ has at least one intersection on \ $[2k\pi-\frac{\pi}{2},2k\pi+\frac{\pi}{2}]$, that is to say, there exists \ $x_0\in [2k\pi-\frac{\pi}{2},2k\pi+\frac{\pi}{2}]$ such that \ $\sin (x_0)=s(x_0)$, so \ $f_1(x_0,\sin (x_0))=0$, and clearly  \ $\cos (x_0)>0$. As \ $f_1(x_0,\sin (x_0))= res(f(x_0,\sin (x_0),c)$,$\sin^2(x_0)+c^2-1,c) =0$ implies \ $f(x_0,\sin (x_0),c)$ has common root with $\sin^2(x_0)+c^2-1$, \ $\sin^2(x_0)+c^2-1=0$ and \ $c>0$ derive \ $c=\cos (x_0)$, i.e. \ $f(x_0,\sin (x_0),\cos (x_0))=0$, \ $x_0$ is a root of \ $F(x)$, which is bigger than the given \ $T$, so the positive roots of \ $F(x)$ is unbounded.

Suppose that \ $c<0$ for \ $x\in(r_0,+\infty)$. As \ $\sin(x)$ changes from \ $1$ to \ $-1$ continuously on \ $[2k\pi+\frac{\pi}{2},2k\pi+\frac{3\pi}{2}]$, then for \ $\forall k\in \mathbb{N}$ such that \ $ 2 k\pi+\pi/2>k$,\ $\sin (x)$ and $s(r)$ have at least one intersection on \ $[2k\pi+\frac{\pi}{2},2k\pi+\frac{3\pi}{2}]$, that is to say, there exists \ $x_0\in [2k\pi+\frac{\pi}{2},2k\pi+\frac{3\pi}{2}]$ such that \ $\sin (x_0)=s(x_0)$, so \ $f_1(x_0,\sin (x_0))=0$, obviously \ $\cos(x_0)<0$. \ $\cos (x_0)$ is the common root of \ $f(x_0,\sin(x_0),c)$ and \ $\sin(x_0)^2+c^2-1$, i.e. \ $f(x_0,\sin (x_0),\cos (x_0))=0$, \ $x_0$ is a root of \ $F(x)$, which is bigger than the given \ $T$, so the positive roots of  \ $F(x)$ is unbounded.

Now suppose that $degree(f,c)$ is even, \ $g(x,s)= f(x,s,c)$ after the substitution \ $c^2=1-s^2$. In this case, \ $F(x)=g(x,sin(x))$ and \ $g\in \mathbb{R}_{alg}[x,y]$ obviously.

Let \ $A$ be a connected branch of $\{x\in \mathbb{R}|g_1(x)= res(s+1,g,s)\times res(s-1,g,s)\times res(g,diff(g,s),s)\neq 0\}$, which is an open interval, \ $(s+1)$ has constant sign on \ $A$ and so does \ $(s-1)$, and for \ $\forall x_1,x_2\in A$, by Lemma 2.1 \ $g(x_1,s)$ and \ $g(x_2,s)$ has the same number of real roots.

If \ $\exists r>r_0$, \ $g_2(s)=g(r,s)$ has a least one real root on \ $(-1,1)$, then \ $\forall r\in (r_0,+\infty)$, \ $g_2(s)=g(r,s)$ has at least one real root on \ $(-1,1)$. From the previous reasoning we can obtain that for \ $k\in \mathbb{N}$ such that \ $2 k\pi-\pi/2>r_0$, $\exists x_0\in [2 k\pi-\pi/2,2 k\pi+\pi/2]$ such that \ $g(x_0,sin(x_0))=0$,  which implies that the positive root of \ $g(x,sin(x))$ is unbounded.

If \ $\exists r>r_0$, \ $g_2(s)=g(r,s)$ has no real root on \ $(-1,1)$, then \ $\forall r\in (r_0,+\infty)$, \ $g_2(s)=g(r,s)$ has no real root on \ $(-1,1)$, so \ $g(r,sin(r))\neq 0$ for \ $\sin(r)\in (-1,1)$. If \ $\sin(r)=1$ or \ $-1$, \ $r=2 k\pi+\pi/2$ or \ $2 k\pi-\pi/2$, as \ $g\in \mathbb{R}_{alg}[x,y]$, we can clare that \ $g(2 k\pi+\pi/2, 1)\neq 0$ and \ $g(2 k\pi-\pi/2, -1)\neq 0$. So $g(r,sin(r))\neq 0$ for \ $\forall r>r_0$, i.e. \ $g(x,sin(x))$ has no real root for \ $x>r_0$, \ $r_0$ is an upper bound of positive roots of \ $F(x)=g(x,sin(x))$.

\end{proof}

To describe the running process of Algorithm 2.1, we present the following examples.

\noindent\textbf{Example 1}
Decide whether the positive roots of  \ $F(x)=\frac{2}{3}x+x\cos(x)-\sin(x)$ is bounded.

Let \ $f(x,s,c)=\frac{2}{3}x+x c-s $, where \ $c$ denotes \ $\cos(x)$ and \ $s$ denotes \ $\sin(x)$,

\ $f_1(x,s)=res(f,s^2+c^2-1,c)=x^2 s^2-\frac{5}{9}x^2-\frac{4}{3}x s+s^2 $,

\ $g(x)=res(f_1, diff(f_1, s), s)=-\frac{20}{9} x^{6}-\frac{56}{9} x^{4}-4 x^{2}$,

\ $realroot(g(x), \frac{1}{10})=[[0,0]]$, which implies that \ $g(x)$ has one root \ $0$.

\ $f_2(s)=subs(x=1, f_1(x))=2 s^{2}-\frac{5}{9}-\frac{4}{3} s$,

\ $realroot(f_2(s), \frac{1}{10})=[[\frac{15}{16}, 1],[-\frac{5}{16},-\frac{1}{4}]]$, \ $f_2$ has \ $2$ real roots, which implies that \ $g_x(s)=f_1(x,s)$ has \ $2$ real roots for \ $\forall x>0$.

As \ $f(x,1,0)=\frac{2}{3} x - 1$,\  $f(x,-1,0)=\frac{2}{3} x + 1$, whose roots are \ $\frac{3}{2}$ and \ $-\frac{3}{2}$, so under the assumption that \ $f(x,s,c)=\frac{2}{3}x+x c-s=0$ and \ $s^2+c^2=1$, if \ $x>\frac{3}{2}$, \ $c<0$ definitely.

Let \ $T$ be a arbitrary given positive real number, \ $T_0=max\{\frac{3}{2},T \}$, let \ $k_0$ be a number such that \ $2k_0\pi + \pi/2>T_0$, then for each natural number \ $k\geq k_0$, \ $\exists x_0\in [2k\pi+\frac{\pi}{2},2k\pi+\frac{3\pi}{2}]$ such that \ $F(x_0)=f(x_0,sin(x_0),cos(x_0))=0$. So \ $x_0$ is a positive root of \ $F(x)$ and \ $x_0>T$, we get that the positive roots of \ $F(x)$ are unbounded.

\noindent\textbf{Example 2}
Decide whether the positive roots of \ $F(x)=\frac{2}{3}x+\frac{1}{3} x\cos(x)-\sin(x)$  is bounded.

Let \ $f(x,s,c)=\frac{2}{3}x+\frac{1}{3} x \cdot c-s $,

%以下作出符号有修订，realroot结果是新做的但之前的也没毛病，对结论没影响
\ $f_1(x,s)=res(f,s^2+c^2-1,c)=\frac{1}{9}x^2 s^2+\frac{1}{3}x^2-\frac{4}{3}x s+s^2 $,

\ $g(x)=res(f_1, diff(f_1, s), s)=\frac{4}{243} x^{6}+\frac{8}{81} x^{4}-\frac{4}{9} x^{2}$,

\ $realroot(g(x), \frac{1}{10})=[[-\frac{111}{64},-\frac{221}{128}],[0,0],[\frac{221}{128},\frac{111}{64}]]$,

\ $f_2(x)=subs(x=2, f_1)=\frac{13}{9} s^{2}+\frac{4}{3}-\frac{8}{3} s$,

\ $realroot(f_2(x), \frac{1}{10})=[]$, which means that\ $f_2(x)$ has no real roots and implies that \ $g_x(s)=f_1(x,s)$ has no real roots for \ $\forall x>\frac{111}{64}$. Then for \ $\forall r>\frac{111}{64}$, \ $sin(r)$ is not the root of  $g_r(s)=f_1(r,s)$,
 which is to say, \ $f_1(r,sin(r))=res(f(r,\sin(r),c),\sin^2(r)+c^2-1,c)\neq 0$. So, \ $\cos(r)$ is not a root of \ $f(r,\sin(r),c)$, otherwise, $\cos(r)$ is the common root of \ $f(r,\sin(r),c)$ and \ $\sin^2(r)+c^2-1$.

So, for \ $\forall r>\frac{111}{64}$,\ $r$ can not be the real roots of \ $F(x)$, i.e. the real roots of \ $F(x)$ are not bigger than \ $\frac{111}{64}$.

\section{Square-free factorization of Mixed Trigonometric polynomial}

Euler Theorem shows that \ $\sin(x)=\frac{e^{I x}-e^{-I x}}{2 I}$, \ $\cos(x)=\frac{e^{I x}+e^{-I x}}{2}$. Hence, an MTP can be expressed as an exponential polynomial in complex field \ $f(x, e^{I x}, e^{-I x})$, where \ $I^2=-1$.

Let \ $y=e^{Ix}$, then an MTP can be reduced to a Laurent polynomial in the form of \ $f[x, y, y^{-1}]$. Denotes \ $\mathbb{LR} := \mathbb{C}[x, y, y^{-1}]$, then for any \ $P\in \mathbb{LR}$, there exists a \ $Q\in \mathbb{C}[x, y]$,\  $P=Q/y^p$, \ $p\in \mathbb{Z}$. If \ $factor(Q)$ is a factorization of \ $Q$, we call \ $factor(Q)/y^p$ or \ $y^{-p} factor(Q)$ a factorization of \ $P$.

By Lindemann Theorem, we have that

\noindent\textbf{Lemma 3.1}
For \ $\forall x\in \mathbb{C}$, there are at least one transcendental number between \ $x, e^{Ix}$.

Adjust Theorem 6 of ref[13] slightly, we get the following lemma and its proof.

\noindent\textbf{Lemma 3.2}
If \ $f_1(x, y), f_2(x, y)\in \mathbb{A}[x, y]$ are co-prime, then \ $F_1(x)=f_1(x, e^{Ix})$ and \ $F_2(x)=f_2(x, e^{Ix})$ have no common roots other than \ $0$.
\begin{proof}
Suppose \ $F_1(x)$ and \ $F_2(x)$ have common root \ $x_0$ and \ $x_0 \neq 0$.

Let \ $g(x)=res(f_1(x,y),f_2(x,y) ,y)$ then \ $g\in \mathbb{A}[x]$ and \ $g(x_0) =0$. As \ $f_1(x, y)$ and \ $f_2(x, y)$ are co-prime, then \ $g(x)$ can not be identically zero, so \ $x_0$ is algebraic. Then \ $h(y)=f_1(x_0,y) \in \mathbb{A}[y]$, so \ $h(e^{Ix_0})=f_1(x_0 ,e^{Ix_0}) =0$, which implies \ $e^{Ix_0}$ is algebraic and contradicts Lemma 3.1.
\end{proof}

\noindent\textbf{Theorem 3.1}
If \ $f(x, y)\in \mathbb{A}[x, y]$ is irreducible, then \ $F(x)= f(x, e^{Ix})$ has no multiple roots other than \ $0$.
\begin{proof}
Let \ $f'(x,y)\in \mathbb{A}[x, y]$ such that \ $F'(x)=f'(x,e^{I x})$, then \ $f'(x,y)=f'_x + I y f’_y$, so \ $degree(f',x)\leq degree(f,x)$, \ $degree(f',y)\leq degree(f,y)$. As \ $f(x,y)$ is irreducible, so \ $f(x,y)$ and \ $f'(x,y)$ are co-prime. Then \ $F(x)$ and \ $F'(x)$ have no common roots other than \ $0$. We conclude that Theorem 3.1 holds.
\end{proof}

\noindent\textbf{Corollary 3.1}
If \ $f(x, y)\in \mathbb{A}[x, y]$ is square-free, then \ $F(x)= f(x, e^{Ix})$ has no multiple roots other than \ $0$.

We extend the operation of complex conjugation to \ $\mathbb{LR}$ as follows.

Given \ $P = \sum \limits_{j=1} ^{n} a_j(x) y^{v_j}\in \mathbb{LR}$, where \ $a_1,\cdots, a_n \in \mathbb{C}[x]$, define its conjugate to be \ $con(P) = \sum^n_{j=1}\overline{a_j(x)} y^{-v_j}$, where  \ $\overline{a_j(x)}$ is the conjugate function of \ $a_j(x)$, and obviously, \ $con(a_j(x))=\overline{a_j(x)}$.

In this section, we extend the coefficients of MTP to complex field. For \ $P = \sum \limits_{j=1} ^{n} a_j(x)y^{v_j} \in \mathbb{LR}$,\ $v\in \mathbb{Q}$, denote \ $LRhom[v](P) =P(x,e^{I v x}) = \sum \limits_{j=1} ^{n} a_j(x)(\cos(v_j v x)+I \sin(v_j v x))$, then \ $LRhom[v](P)$ is a mixed trigonometric-polynomial with coefficients in complex field. \ $LRhom[1](P)$ is abbreviated as \ $LRhom(P)$.

For an MTP \ $F(x)=f(x,\sin(x),\cos(x))$, let \ $P(x,y)=f(x,\frac{y-1/y}{2 I},\frac{y+1/y}{2})\in \mathbb{LR}$, then \ $LRhom(P)=F(x)$ obviously.

\noindent\textbf{Lemma 3.3} For \ $P = \sum \limits_{j=1} ^{n} a_j(x) y^{v_j}\in \mathbb{LR}$, \ $LRhom[v](con(P))=\overline{LRhom[v](P)}$ for $\forall v\in \mathbb{Q}$.
\begin{proof}
$LRhom[v](con(P))= LRhom[v](\sum^n_{j=1}\overline{a_j(x)} y^{-v_j})=  \sum \limits_{j=1} ^{n} \overline{a_j(x)}(\cos(-v_j v x)+I \sin(-v_j v x)) = \sum \limits_{j=1} ^{n} \overline{a_j(x)}(\cos(v_j v x)-I \sin(v_j v x))= \sum \limits_{j=1} ^{n} \overline{a_j(x) (\cos(v_j v x)+I \sin(v_j v x))}= \overline{LRhom[v](P)}$
\end{proof}

By Lemma 3.3 we get that

\noindent\textbf{Lemma 3.4}
If \ $P=con(P)$, then \ $LRhom[v](P)$ is real-valued for $\forall v\in \mathbb{Q}$.

\noindent\textbf{Theorem 3.2} % 2.4 for reachibility
If \ $P = \sum_{j=1}^{n}a_j(x) y^{v_{j}}\in \mathbb{LR}$ such that \ $P=con(P)$, \ $P$ can be factorized as \ $C y^p P_1^{r_1} \cdots P_n^{r_n}$, where \ $C\in \mathbb{C}$ and \ $P_1,\dots,P_n\in \mathbb{C}[x,y]$ are square-free and pairwisely co-prime, then for each \ $i$,

1) there exist \ $C_i\in \mathbb{C}$,\ $p_i \in \mathbb{Z}$, such that \ $P_i=C_i y^{p_i} con(P_i)$;

2) $f_i(x)=LRhom(P_i Z^{-p_i/2} C_i^{-1/2})$ is real-valued or pure imaginary, where \ $f_i$ has no multiple root other than \ $0$ for \ $i=1,\dots, n$,\ $f_i$ and \ $f_j$ have no common root other than \ $0$ for \ $i\neq j$;
\
3) \ $LRhom(P) = C_0 f_1^{r_1} \cdots f_n^{r_n}$, where \ $C_0 = C (C_1)^{r_1/2} \dots (C_n)^{r_n/2}$.

\begin{proof}
1)\ For \ $n=1$, $C Z^p (P_1)^{r_1}= \overline{C} Z^{-p} con(P_1)^{r_1}$, thus \ $P_1 = (\overline{C}/C)^{1/r_1} Z^{-2 p/r_1} con(P_1)$, the conclusion holds.

Suppose the conclusion holds for \ $n=m-1$, that is to say, there exist \ $C_i, p_i$ such that \ $P_i=C_i Z^{p_i} con(P_i)$ for \ $1<i<m$.  Let \ $n=m$,

\ $P=C Z^p (P_1)^{r_1} \cdots (P_{m-1})^{r_{m-1}} P_m^{r_m}=C Z^p (C_1)^{r_1} Z^{p_1 r_1} con(P_1)^{r_1} \cdots (C_{m-1})^{r_{m-1}} Z^{p_{m-1} r_{m-1}} con(P_{m-1})^{r_{m-1}} (P_m)^{r_m}$

\ $=C_0 Z^{p_0} con(P_1)^{r_1} \dots con(P_{m-1})^{r_{m-1}} P_m^{r_m}$, where \ $C_0=C (C_1)^{r_1} \dots (C_{m-1})^{r_{m-1}}$, \ $p_0=p+p_1 r_1+\cdots+p_{m-1} r_{m-1}$.

On the other hand, \ $con(P)=\overline{C} Z^{-p} con(P_1)^{r_1} \cdots con(P_{m-1})^{r_{m-1}} con(P_m)^{r_m}$.

So, by the assumption \ $P = con(P)$, we get that if \ $con(P_1)\neq 0, \cdots,con(P_{m-1})\neq 0$, \ $P_m=C_m y^{p_m} con(P_m)$, where \ $C_m=(\overline{C}/C_0)^{1/r_m}$,\ $p_m=(-p-p_0)/r_m$. As \ $con(P_1),\cdots,con(P_{m-1})$ have finite zeros at most, \ $P_m$ and \ $con(P_m)$ are both continuous, \ so $P_m=C_m y^{p_m} con(P_m)$ holds.

2)\ $P_i=C_i Z^{p_i} con(P_i)$ implies \ $P_i^2=C_i Z^{p_i} con(P_i) P_i$, so \ $P_i^2 C_i^{-1} Z^{-p_i}=con(P_i) P_i$. Let \ $Q_i=P_i Z^{-p_i/2} C_i^{-1/2}$, then \ $Q_i^2=con(P_i) P_i$, so that $Q_i^2=con(Q_i^2)$.

Let \ $f_i(x) = LRhom[u_1, \dots, u_r; v_1,\dots,v_s](Q_i)$, then \ $f_i(x)^2$ is real-valued due to Lemma 3.4, i.e. \ $f_i(x)$ is real-valued or pure imaginary.

It is clearly that \ $f_i$ has no multiple root other than \ $0$ by Theorem 3.1,\ $f_i$ and \ $f_j$ have no common root other than \ $0$ for \ $i\neq j$ by Lemma 3.2.

3)\ As \ $P=con(P)$, then \ $degree(P, y)=degree(P,y^{-1})$,  denoted by \ $q$. It is trivial that \ $p=-q$ and  \ $degree(Q = P_1^{r_1} \cdots P_n^{r_n},y)=2 q$. Since \ $P_i=C_i y^{p_i} con(P_i)$, \ $p_i = degree(P_i,y)=degree(con(P_i),y^{-1})$, hence \ $p_1+\dots+p_n = degree(Q,y)$, i.e. \ $p= -q=-(p_1+\dots+p_n)/2$.

So, \ $ C_0 f_1^{r_1} \cdots f_n^{r_n} = C (C_1)^{r_1/2} \dots (C_n)^{r_n/2} (LRhom(P_1 y^{-p_1/2} C_1^{-1/2}))^{r_1} \dots (LRhom(P_n y^{-p_n/2} C_n^{-1/2}))^{r_n}$

\ \ $=  C (LRhom[u_1, \dots, u_r; v_1,\dots,v_s](P_1 y^{-\frac{p_1}{2}} ))^{r_1} \dots (LRhom[u_1, \dots, u_r; v_1,\dots,v_s](P_n y^{-\frac{p_n}{2}}))^{r_n}$

\ \ $= C LRhom[u_1, \dots, u_r; v_1,\dots,v_s](y^{-\frac{p_1}{2} r_1 - \dots -\frac{p_n}{2} r_n} P_1^{r_1} \dots P_n ^{r_n})$

\ \ $=LRhom[u_1, \dots, u_r; v_1,\dots,v_s](C y^p P_1^{r_1} \dots P_n ^{r_n}) = LRhom[u_1, \dots, u_r; v_1,\dots,v_s](P)$
\end{proof}

\noindent\textbf{Corollary 3.2 2.4 for reachability}
For each trigonometric-exponential polynomial \ $F(x)=f(x,e^{{u_1} x}, \dots, e^{{u_r} x}, \sin(v_1 x),\dots, \sin(v_s x), \cos(v_1 x),\dots, \cos(v_s x))$, where \ $f\in \mathbb{R}_{alg}[x, y_1, \dots, y_r, z_1, \dots, z_s, w_1, \dots, w_s]$ there exists \ $C\in \mathbb{R}$ and real-valued  trigonometric-exponential polynomials \ $\{f_i\}$ such that \ $f_i$ has no multiple root other than \ $0$ ,\ $f_i$ and \ $f_j$ have no common root other than \ $0$ for \ $i\neq j$, \ $F(x)= C f_1(x)^{r_1} \dots f_n(x)^{r_n}$.
\begin{proof}

  Making substitution for F(x), \ $e^{{u_i} x}= y_i$ for \ $i=1, \dots, r$, \ $\sin(v_i x)=\ frac{z_i - 1/z_i}{2 I}$, \ $\cos(v_i x)=\ frac{z_i + 1/z_i}{2}$ for \ $i=1, \dots, s$, yields \ $P \in \mathbb{LR}$ such that \ $LRrom[u_1, \dots, u_r; v_1,\dots,v_s](P)=F(x)$.

  By Theorem 2.4, there exists \ $C\in \mathbb{C}$ and  trigonometric-exponential polynomials \ $\{f_i\}$ such that \ $f_i$ has no multiple root other than \ $0$ ,\ $f_i$ and \ $f_j$ have no common root other than \ $0$ for \ $i\neq j$, \ $F(x)= C f_1(x)^{r_1} \dots f_n(x)^{r_n}$, each \ $f_i$ is real-valued or pure imaginary.

  If all \ $f_i$s are real-valued, then Corollary 2.4 holds. If \ $f_i$ is pure imaginary, let \ $f_i'=f_i/I$ and \ $C'=C I^{r_i}$, then  $f_i'$ is real-valued and \ $F(x)= C' f_1(x)^{r_1} \dots {f_i'(x)}^{r_i} \dots f_n(x)^{r_n}$.

  Now, \ $F(x)$ and all \ $f_i(x) (or f_i'(x))$s are real-valued, so the constant \ $C(or C')$ must be a real number. That is to say the corollary holds.

  \end{proof}

\noindent\textbf{Algorithm 3.1}Square-free factorization of MTP

Input:\ an MTP \ $F(x) =f(x,\sin(x),\cos(x)), f\in \mathbb{R}_{alg}[x,y,z]$;

Output:\ $F(x) = C_0 F_1(x)^{n_1} \cdots F_m(x)^{n_m}$; where\ $C_0 \in \mathbb{R}$ and \ $F_i(x)$ is a real-valued MTP and which has no multiple roots other than \ $0$ for \ $i = 1,\cdots, m$ , \ $F_i(x)$ and \ $F_j(x$) have no common roots other than \ $0$ for $i \neq j$.

1)\ $P\gets f(x,\frac{y-1/y}{2 I},\frac{y+1/y}{2})$;

2)\ $P\_lst \gets factor(P)=C y^p P_1^{r_1} \cdots P_n^{r_n}$; where \ $C\in \mathbb{C}$, \ $P_1,\cdots,P_n\in \mathbb{C}[x,y]$ are square-free and pairwisely co-prime.

3)\ $F\gets 1$; \ $C_0\gets C$;

4)\  for \ $i$ form \ $1$ to \ $m$

\ \ \ \ 4.1)\ $g\gets P_i$;

\ \ \ 4.2)\ $h\gets g/con(g)$; where \ $h$ is of form \ $C_i y^{p_i}$, \ $LRhom(g^2/h)=LRhom(g\times  con(g))$ is real-valued;

\ \ \ \ 4.3)\ $f_i\gets LRhom(g/(C_i^{1/2} y^{p_i/2}))$; where \ $f_i$ is real or pure imaginary;

\ \ \ \ 4.4)\ if \ $f_i$ is pure imaginary, then \ $f_i\gets f_i/I$,\ $C_0\gets C_0 I^{r_i}$;

\ \ \ \ 4.5)\ $F\gets F f_i^{r_i}$;

\ \ \ \ 4.6)\ $C_0\gets C_0 (C_i^{1/2})^{r_i}$;

5)\ return \ $C_0 F$.

\noindent\textbf{Example 3}
Decide whether \ $f(x)=\frac{2}{3} x+x \cos(x)-\sin(x)$ has multiple roots.

Let \ $\cos(x)=\frac{e^{I x}+e^{-I x}}{2}$,\ $\sin(x)=\frac{e^{I x}-e^{-I x}}{2 I}$,\ $y=e^{I x}$,

then \ $f=\frac{2 x}{3}+x \frac{y+1/y}{2}-\frac{y-1/y}{2 I}$,

\ $factor(f)=\frac{1}{6} \frac{4 x y+3 x y^{2}+3 x+3 I y^{2}-3 I}{y}$, \ $4 x y+3 x y^{2}+3 x+3 I y^{2}-3 I$ is irreducible, which implies that \ $f(x)$ has no multiple roots.

\noindent\textbf{Example 4}
Decide whether \ $f(x)=1-\sin^3(x)$ has multiple roots, and if so, do factorization.

Let\ $\sin(x)=\frac{e^{I x}-e^{-I x}}{2 I}$, \ $y=e^{I x}$.

Then $f(x)=g(y)=1-\frac{1}{8}I(y-\frac{1}{y})^3= \frac{-\frac{1}{8}I(y^4+2Iy^3-6y^2-2Iy+1)(y-I)^2}{y^3}=C P_1 (P_2)^2$, where \ $C=-\frac{I}{8}$, \ $P_1 = y^4+2 I y^3-6 y^2-2 I y+1$, and \ $P_2=y-I$.

 We get that \ $y=I$ is a multiple root of \ $g(x)$, and then \ $ \{2 k \pi+\frac{\pi}{2}, k\in \mathbb{Z} \}$ are multiple roots of \ $f(x)$.

As \ $con(P_1)=\frac{1}{y^4}-\frac{2 I}{y^3}-\frac{6}{y^2}+\frac{2I}{y}+1$, so \ $P_1=c_1 y^4 con(P_1)$, where \ $c_1=1$.

Let
\begin{equation}
\begin{aligned}
g_1 &= LRhom(P_1/(y^4)^{1/2}) = LRhom((y^4+2 I y^3-6y^2-2 I y+1)/y^2) \\
&=LRhom(y^2+2 I y-6-2 I y^{-1}+y^{-2})\\ &=e^{2 I x}+2 e^{I \frac{\pi}{2}} e^{I x} -6-2 e^{I{\frac{\pi}{2}}} e^{-I x}+e^{-2 I x} \\
&=\cos(2x)+I \sin(2x)+2 \cos(\frac{\pi}{2}+x)+2 I \sin(\frac{\pi}{2}+x)-6-2\cos(\frac{\pi}{2}-x)\\
&-2 I \sin(\frac{\pi}{2}-x) +\cos(-2x)+I \sin(-2x)\\
&=\cos(2x)+I \sin (2x)-2\sin(x)+2 I \cos(x)-6-2 \sin(x)-2 I \cos(x)+\cos(2x)-I \sin(2x)\\
&= 2 \cos(2x)-4 \sin(x)-6
\end{aligned}
\end{equation}

As \ $con(P_2) = \frac{1}{y} + I$, then \ $P_2 = c_2 y con(P_2)$, where \ $c_2 = -I = e^{-I \frac{\pi}{2}}$.

%g2已作修订，请复核
Let
\begin{equation}
\begin{aligned}
g_2 &= LRhom(P_2/(-I y)^{\frac{1}{2}} )= LRhom((y - I)/(-I y)^{\frac{1}{2}} ) \\
&= (e^{I x} - e^{I \frac{\pi}{2}})/e^{-I \frac{\pi}{4}+I \frac{x}{2}}\\
&= e^{I \frac{x}{2}+I \frac{\pi}{4}} - e^{I \frac{3\pi}{4}-I \frac{x}{2}}\\
&= \cos(\frac{x}{2} + \frac{\pi}{4}) + I \sin(\frac{x}{2} + \frac{\pi}{4}) - (\cos(\frac{3\pi}{4} - \frac{x}{2}) + I \sin(\frac{3\pi}{4} -\frac{x}{2}))\\
&= \cos(\frac{x}{2} + \frac{\pi}{4}) + I \sin(\frac{x}{2} + \frac{\pi}{4}) + \cos(\frac{x}{2}+\frac{\pi}{4}) - I \sin(\frac{x}{2}+\frac{\pi}{4})\\
&= 2 \cos(\frac{\pi}{4}+\frac{x}{2})
\end{aligned}
\end{equation}

Let \ $c_0=C c_1^{1/2} (c_2^{1/2})^2 = -\frac{I}{8}(-I)=-1/8$.

Then \ $f(x)=c_0 g_1 (g_2)^2 =-\frac{1}{8}(2 \cos(2 x)-4 \sin(x)-6) \cos^2(\frac{\pi}{4}+\frac{x}{2})=-\frac{1}{2} (\cos(2x)-2\sin(x)-3) (\cos(\frac{x}{2}) - \sin(\frac{x}{2}))^2=-\frac{1}{2} f_1 (f_2)^2$, where \ $f_1=\cos(2x)-2\sin(x)-3, f_2=\cos(\frac{x}{2}) - \sin(\frac{x}{2})$.

\section{Sign-deciding of Mixed Trigonometric Polynomial}
If the positive roots of a square-free MTP are unbounded, it has no constant sign on $(0,+\infty)$ definitely.
So we need only to discuss the sign-deciding of square-free MTP on a bounded interval.

Ref[8] has designed a complete algorithm \ $Deciding\_arctan\_polynomial$ for deciding the sign of  inverse tangent function polynomial $f(x,arctan(x))$ on \ $(0, 1)$. Each mixed trigonometric-polynomial can be transformed to a tangent function polynomial \ $f(x,\tan(\frac{x}{2}))$. Let \ $y = \tan(\frac{x}{2})$, then \ $x = 2 \arctan(y)$ and \ $f(x,\tan(\frac{x}{2})) = f(2 \arctan(y), y)$, thus tangent function polynomial is reduced to an inverse tangent function polynomial. So it is trivial to decide the sign of an MTP on \ $(0,\frac{\pi}{2})$.(In general, we don't care whether the right end of the interval is open or closed.)

For \ $T>\frac{\pi}{2}$, to decide the sign of MTP on (0,T), the tan-half angle substitutions \ $\cos(2 t)=1-2\sin^2 (t)$ and $\sin(2 t)=2\sin (t) \cos (t)$ can reduce \ $f(x)$ to another MTP \ $g(t)$ such that \ $0<t<\frac{\pi}{2}$.

\noindent\textbf{Algorithm 4.1} \ $Decide\_trigonometric\_polynomial$

Input:\ $F(x)=f(x,\sin(x),\cos(x))$, where $f\in \mathbb{R}_{alg}[x,y,z]$ and is square-free, a constant \ $T$;

Output:\ The sign of \ $F(x)$ on \ $(0,T)$.

1)\	$p \gets [T/(\frac{\pi}{2})]+1$; where \ $[x]$ denotes the maximum integer not bigger than \ $x$

2)\	$G(t) \gets F(p t)$, and transform \ $G(t)$ to the form \ $g(t,\sin(t),\cos(t))$;

3)\	$G(t)\gets g(t, \frac{2\tan(\frac{t}{2})}{1+\tan^2(\frac{t}{2})},\frac{1-\tan^2(\frac{t}{2})}{1+\tan^2(\frac{t}{2})})$;

4)\	$H(y)\gets subs(\tan(\frac{t}{2})=y, t=2\arctan(y),G(t))$;

5)\	Return \ $Deciding\_arctan\_polynomial(H(y))$.

For a general MTP on $(0,+\infty)$, we design the following algorithm.

\noindent\textbf{Algorithm 4.2}\ $Decide\_general\_trigonometric\_polynomial$

Input:\ $F(x)=f(x,\sin(x),\cos(x))$, \ $f\in \mathbb{R}_{alg}[x,y,z]$;

Output:\ The sign of \ $F(x)$ on \ $(0,+\infty)$;

1)\ Run Algorithm 3.1 to factorize \ $F(x)$, \ $F(x) = c f_1^{d_1} f_2^{d_2} \cdots f_n^{d_n}$; Where \ $c$ is a real constant, \ $f_1, f_2,\cdots, f_n$ are real-valued MTP and has no multiple roots other than \ $0$ for \ $i = 1,\cdots, n$ , \ $f_i(x)$ and \ $f_j(x$) have no common roots other than \ $0$ for $i \neq j$.

2)\ $p \gets 0$, $q \gets 0$; Where \ $p$ denotes the number of negative terms, \ $q > 0$ shows that the inequality is not strict.

3)\  $m \gets 1$.

4)\ if \ $d_m$ is odd, then

\ \ \ \ Run Algorithm 2.1 to decide if the positive roots of \ $f_m$ are bounded,

\ \ \ \ 4.1)\ If unbounded, then return \ $0$, the algorithm ends; where \ $0$ means that \ $F(x)$ has no constant sign on \ $(0,+\infty)$.

\ \ \ \ 4.2)\ If bounded and \ $T$ is a bound, then \ $s \gets Decide\_trigonometric\_polynomial(f_m,T)$

\ \ \ \ \ \ If \ $s=-1$ then \ $p \gets p + 1$; where \ $s=-1$ implies that \ $f_m<0$ holds on \ $(0,T)$.

\ \ \ \ \ \ If \ $s=0$ then return \ $0$,the algorithm ends; where \ $s=0$ implies that \ $f_m$ has no constant sign \ $(0,T)$.

\ \ \ elif \ $d_m$ is even, then

\ \ \ \ 4.3)\ If the positive roots of \ $f_m$ are unbounded, then \ $q \gets q + 1$;

\ \ \ \ 4.4)\ The positive roots of \ $f_m$ are bounded and \ $T$ is a bound, then \ $s \gets Decide\_trigonometric\_polynomial(f_i,T)$,

\ \ \ \ \ \ If \ $s=0$ then \ $q \gets q + 1$.

5)\  $m \gets m + 1$, if \ $m < n$ then goto 4)

6)\ 6.1)\ If \ $c > 0$, \ $p$ is even and \ $q = 0$, then inequality \ $F(x) > 0$ holds;

\ \ \ \ 6.2)\ If \ $c > 0$, \ $p$ is even and \ $q > 0$, then inequality \ $F(x) \geq 0$ holds;

\ \ \ \ 6.3)\ If \ $c > 0$, \ $p$ is odd and \ $q > 0$, then inequality \ $F(x) < 0$ holds;

\ \ \ \ 6.4)\ If \ $c > 0$, \ $p$ is odd and \ $q > 0$, then inequality \ $F (x) \leq 0$ holds;

For \ $c < 0$, the conclusions can be drawn in the same manner.

\noindent\textbf{Example 5}
Decide the sign of \ $f(x)=\frac{2}{3}x+x\cos(x)-\sin(x)$ on \ $(0,+\infty)$.

By Example 1, the positive roots of \ $f(x)$ are unbounded. By Example 3, \ $f(x)$ has no multiple root. So \ $f(x)$ has no constant sign on \ $(0,+\infty)$.

\noindent\textbf{Example 6}
Decide the sign of \ $f(x)=\frac{2}{3}x+\frac{1}{3} x\cos(x)-\sin(x)$ on \ $(0,+\infty)$.

%依据新的例2 \frac{111}{64}
By Example 2, the positive roots of \ $f(x)$ are less than \ $\frac{111}{64}$, so we need only to discuss the problem on \ $(0,\frac{111}{64})$.

As \ $\frac{111}{64}>\frac{\pi}{2}$, let \ $g(y)=f(2y)=\frac{4}{3}y + \frac{2}{3}y(1-2 \sin(y)^2) - 2\sin(y)\cos(y)$, where \ $0<y<\frac{111}{128}<\frac{\pi}{2}$.

Let \ $t=\tan(\frac{y}{2})$, then \ $0<t<1$, \ $\sin(y)=\frac{2 t}{1+t^2}$, \ $\cos(y)=\frac{1-t^2}{1+t^2}$, hence \ $g(y)= \frac{2}{3} \frac{3 y-2y t^2+3 y t^4 -6 t +6 t^3}{(1+t^2)^2} $.

Now it is necessary only to determine the sign of \ $g_1(y,t)= 3 y-2y t^2+3 y t^4 -6 t +6 t^3$. Obviously \ $y=2\arctan(t)$, so \ $g_1(y,t)= 3 \arctan(t)-2\arctan(t) t^2+3 \arctan(t) t^4 -6 t +6 t^3 $, denoted as \ $G(t)$, the algorithm \ $Deciding\_arctan\_polynomial$ of Ref[8] can declare that \ $G(t)>0$ holds on \ $(0,1)$, and so \ $g(y)=f(2y)>0$ on \ $(0,\frac{\pi}{2})$, \ $f(x)>0$ on \ $(0,\pi)$. As \ $f(x)$ has constant sign on \ $(\frac{111}{64},+\infty)$, we get that \ $f(x)>0$ holds on \ $(0,+\infty)$.

\section{Conclusion}
In this paper, we present algorithms to determine the boundedness of MTP's positive roots and to decompose MTPs without multiple roots, and then propose a procedure to decide the sign of MTP on $(0,+\infty)$. Furthermore, the algorithms presented are much useful for similar problems.

\bibliographystyle{elsarticle-harv}

\begin{thebibliography}{30}
\bibitem[]{}
Bercu, G., 2016. Pade approximant related to remarkable inequalities involving trigonometric functions. J. Inequal. Appl. 99.
\bibitem[]{}
Bercu, G., 2017. The natural approach of trigonometric inequalities-Pade approximant. J. Math. Inequal. 11 (1), 181-191.
\bibitem[]{}
Chen, C-P., 2012. Sharp Wilker and Huygens type inequalities for inverse trigonometric and inverse hyperbolic functions. Integral Transforms Spec. Funct. 23 (12), 865-873.
\bibitem[]{}
B. Male\v{s}evi\'c, M. Makragi\'c: A Method for Proving Some Inequalities on Mixed Trigonometric Polynomial Functions, Journal of Mathematical Inequalities, Volume 10, Number 3 (2016), 849-876 (DOI: 10.7153/jmi-10-69).
\bibitem[]{}
B. Male\v{s}evi\'c, B. Banjac, I. Jovovi\'c: A proof of two conjectures of Chao-Ping Chen for inverse trigonometric functions, Journal of Mathematical Inequalities, Volume 11, Number 1 (2017), 151-162 (DOI: 10.7153/jmi-11-15).
\bibitem[]{}
B. Male\v{s}evi\'c, T. Lutovac, B. Banja\'c: A proof of an open problem of Yusuke Nishizawa for a power-exponential function, Journal of Mathematical Inequalities, Volume 12, Number 2 (2018), 473-485 (DOI: 10.7153/jmi-2018-12-35).
\bibitem[]{}
Mortici, C., 2011. The natural approach of Wilker-Cusa-Huygens inequalities. Math. Inequal. Appl. 14 (3), 535-541.
\bibitem[]{}
Chen shiping, Liu zhong,Automated Proof of Mixed Trigonometric-polynomial Inequalities,Journal of Symbolic Computation. 101C (2020) pp. 318-329
\bibitem[]{}
Chonev, V., Ouaknine, J., Worrell, J.,  On the Skolem problem for continuous linear dynamical systems.2016. In: Chatzigiannakis
\bibitem[]{}
Jamie Mulholland and Michael Monagan. 2001. Algorithms for trigonometric polynomials. In Proceedings of the 2001 international symposium on Symbolic and algebraic computation(ISSAC'01). Association for Computing Machinery, New York, NY, USA, 245-252. DOI:https://doi.org/10.1145/384101.384135
\bibitem[]{}
Achim Schweikard ,Trigonometric polynomials with simple roots,Information Processing Letters 39 (1991) 231-236
\bibitem[]{}
Yang, Lu, Xia, Bican, 2008. Inequality Automated Proving and Discovery. Science Press.
\bibitem[]{}
McCallum, S., Weispfenning, V., 2012. Deciding polynomial-transcendental problems. J. Symb. Comput. 47, 16-31.

\end{thebibliography}

\end{document}